\documentclass[11pt]{amsart}
\usepackage{amsmath}

\usepackage{amssymb}
\usepackage{latexsym}
\usepackage{amscd}
\usepackage{enumerate}
\usepackage[all]{xy}
\input xy
\xyoption{all}
\usepackage[bookmarks=false]{hyperref}

\numberwithin{equation}{section}

\def\proof{\noindent\emph{Proof:}\hs}

\newcommand{\om}{\omega}

\newcommand{\tens}{\otimes}

\def\ov{\overline}
\def\Z{\mathbb{Z}}
\def\Q{\mathbb{Q}}
\def\R{\mathbb{R}}
\def\C{\mathbb{C}}
\def\G{PSL_2(\C)}

\def\H{\mathcal{H}}

\def\H{\mathcal{H}}

\def\C{\mathbb{C}}

\def\M{\mathcal{M}}

\def\hra{\hookrightarrow}

\def\la{\langle}
\def\wt{\widetilde}
\def\ra{\rangle}

\def\({\left(}
\def\){\right)}
\def\/{{/\!\!/}}
\def\sra{\rightarrow}

\def\CC{{\mathcal{C}}}

\def\MMC{{\mathcal{M}}}
\def\coker{{\text{coker }}}

\def\delbar{\ov{\partial}}

\def\hsp{\hspace{0.2cm}}
\def\hspa{\hspace{0.3cm}}

\def\and{\hsp\text{and}\hsp}
\def\ov{\overline}

\def\hra{\hookrightarrow}

\def\and{\hspa\text{and}\hspa}

\def\proof{\noindent\emph{Proof:}\,\,\,\,\,}
\def\bs{\backslash}

\newtheorem{defn}{Definition}[section]
\newtheorem{theorem}{Theorem}[section]

\newtheorem*{theorem*}{Theorem}

\newtheorem{cor}[theorem]{Corollary}

\newtheorem{lem}[theorem]{Lemma}
\newtheorem{prop}[theorem]{Proposition}
\theoremstyle{definition}

\def\begineq{\begin{equation}}
\def\endeq{\end{equation}}
\def\begineqno{\begin{equation*}}
\def\endeqno{\end{equation*}}
\def\beginit{\begin{itemize}}
\def\endit{\end{itemize}}
\def\begineqn{\begin{eqnarray}}
\def\endeqn{\end{eqnarray}}

\begin{document}

\title{The Weinstein Conjecture for Hamiltonian fibrations.}

\author{Cl\'ement Hyvrier}
\address{
Universit\'e de Montr\'eal\\
D\'epartement de Math\'ematiques et de Statistiques\\
 Pavillon Andr\'e-Aisenstadt\\
  2920, chemin de la Tour\\
 Montr\'eal (Qu\'ebec) H3T 1J4\\
Canada
}
\email{
hyvrier@dms.umontreal.ca
}

\begin{abstract} In this note we extend to non trivial  Hamiltonian fibrations over symplectically uniruled manifolds a result of  Lu's, \cite{Lu},  stating that any trivial symplectic product of two closed symplectic manifolds with one of them being  symplectically uniruled   verifies the Weinstein Conjecture for  closed separating hypersurfaces of contact type, under certain technical conditions.  The proof is based on   the product formula for Gromov-Witten invariants ($GW$-invariant)  of Hamiltonian fibrations derived in \cite{H}.
\end{abstract}

\maketitle

\section{Introduction} This note is based on the following observation regarding  the existence of a non-vanishing genus zero Gromov-Witten invariant  with a least one point as constraint, namely \emph{symplectic uniruledness} (see Definition \ref{sympunirul}), of closed connected Hamiltonian fibrations over symplectic bases,
$$(F^{2n_F},\om)\stackrel{\iota}{\hra} P^{2n_P}\stackrel{\pi}{\sra} (B^{2n_B},\om_B),$$  
which are  \emph{rationally cohomologically split}, i.e. where the map $\iota$ induces an injective map in rational cohomology.
 The observation can be stated as follows:
\vspace{0.3cm}
\begin{itemize}\item[\textbf{(O)}:] \emph{For   $c$-splitting Hamiltonian fibrations,  symplectic uniruledness of the base is sufficient to ensure  symplectic uniruledness of the total space.}
\end{itemize} 
\vspace{0.3cm}
 
 The $c$-splitting hypothesis actually holds in many cases \cite{LM},\cite{Bl},\cite{Kedra},\cite{H}. It is in fact  conjectured by McDuff and Lalond, \cite{LM}, that  every Hamiltonian fibration is rationally cohomology split. 

Although we think \textbf{(O)}  should hold in full generality, we  only provide a proof under the following  technical assumptions.  First, we assume that $(F,\om)$ is  \emph{semi-positive relative to the total space}, i.e. $(F,\om)$ verifies condition  \eqref{ssp} described in  \textsection 2.3.2. This condition, verified by  Hamiltonian fibrations of real dimension at most six, is a fibered analog of the "standard" semi-positivity  condition used in \cite{MS}.

 Secondly, we assume that the base is symplectically uniruled for a class  $\sigma_B\in H_2(B;\Z)$ which \emph{only admits simple decompositions for some $\om_B$-tame almost complex structure $J_B$ on $B$}.  This condition, explicited in Definition \ref{irred}, guaranties that, for suitable almost complex structure $J_P$ on $P$ lifting $J_B$,  there are no simple stable $J_P$-holomorphic maps in $P$ projecting to a non-simple stable $J_B$-holomorphic maps in $B$ representing $\sigma_B$.  This is in particular realized by primitive classes admitting pseudo-holomorphic representatives. Here is our main result:

\begin{theorem}\label{mainobservation} Let $(F,\om)\stackrel{\iota}{\hra} P\stackrel{\pi}{\sra} (B,\om_B)$ be a $c$-splitting Hamiltonian fibration. Assume $(F,\om)$ is semi-positive relatively to $P$  and that $(B,\om_B)$ is symplectically uniruled for some class $\sigma_B\in H_2(B;\Z)$ admitting only simple decompositions, then $P$ is also symplectically uniruled. 
\end{theorem}

This generalizes a result in  \cite{H} stating  that symplectic uniruling holds for total spaces of Hamiltonian fibrations over symplectically rationally connected base (for a class $\sigma_B$ as above), i.e. when the base has a non-vanishing genus zero Gromov-Witten invariant  with a least two points as constraints. However, in that latter case, no $c$-splitting assumption is required.

 One of the main ingredients of the proof of Theorem \ref{mainobservation} is to use \emph{fibered almost complex structure on $P$} in order to relate pseudo-holomorphic maps in $P$ with pseudo-holomorphic maps in $B$. Roughly speaking, a fibered almost complex structure $J_P$ on $P$ is an almost complex structure lifting some $\om_B$-tame almost complex structure $J_B$ on $B$, and preserving both horizontal and vertical directions in some  splitting of $TP$ induced by a Hamiltonian connection on $P$ (cf Definition \ref{fiberedacs}). This geometrical setup enables to relate, generically,  some  $GW$-invariant of $P$ with some  $GW$-invariant of $B$ via a product-type formula involving $GW$-invariants of some Hamiltonian fibration over $S^2$. The technical assumptions of the theorem are there to ensure that transversality can be realized within the realm of fibered almost complex structures.
 
In a sense,  \textbf{(O)} is complementary to the following result of Ruan and Li:

\begin{theorem}\label{SympDiv}(\cite{RuanLi}, Proposition 2.10) If  $(F,\om)\stackrel{\iota}{\hra}P\stackrel{\pi}{\sra} B$ is a $c$-splitting Hamiltonian fibration, then $P$ is symplectically uniruled for a class $\sigma\in H_2(P;\Z)$ such that $\pi_*\sigma=0$ if and only if $(F,\om)$ is symplectically uniruled.
\end{theorem}

 Note that the theorem above holds  in great generality.  It is believed that the \emph{ad hoc} hypothesis of Theorem \ref{mainobservation} can be ruled out using virtual perturbations (see \cite{CL}, \cite{LT}, \cite{Rvirt}, amongst others). The problem is to show that one can make perturbations compatibly with the fiber bundle projection $\pi$.  Removing these assumptions is  part of a joint  work in progress  with Shengda Hu.      
 
In some circumstances, we can remove the assumption   on $\sigma_B$.  For instance, it  can be removed if the base is a smooth projective manifold $(B,J_B,\om_B)$. In such case,  Ruan and Koll\`{a}r showed that \emph{strong symplectic uniruledness} is equivalent  to \emph{projective uniruledness}, i.e. that there is a rational curve passing through every point of $B$. The following proposition then follows from Theorem \ref{mainobservation}.

 \begin{cor} \label{cormainobs} Let $(F,\om)\stackrel{\iota}{\hra} P\stackrel{\pi}{\sra} B$ be a $c$-splitting Hamiltonian fibration over  a smooth projective manifold $(B,J_B,\om_B)$.   Assume $(F,\om)$ is semi-positive relatively to $P$ and  that  $(B,J_B,\om_B)$ is uniruled, then $P$ is strongly symplectically uniruled. 
\end{cor}


In view of Theorem \ref{mainobservation} and Theorem \ref{SympDiv}, it is natural to ask the following question: "If neither $F$ or $B$ is symplectically uniruled, can we conclude that $P$ is not symplectically uniruled ?". We give a partial answer:
 
 \begin{theorem} \label{nouniruling} Assume $(B,J_B,\om_B)$ is a smooth  projective manifold and that $P$ c-splits. Suppose that  neither $(B,J_B,\om_B)$, nor $(F,\om)$, is symplectically uniruled, then $P$ is not symplectically uniruled for any (generic) fibered almost complex structure $J_P$ on $P$ lifting $J_B$.
\end{theorem} 
  
 \proof 
Assume by contradiction that $P$ is symplectically uniruled for some class $\sigma\in H_2(P;\Z)$. Two possibilities may occur:  $\pi_*\sigma$ is zero or not. In case $\pi_*(\sigma)=0$, we can use the c-splitting hypothesis to apply Theorem \ref{SympDiv} and obtain a contradiction. Now, let's consider the case where $\pi_*\sigma\neq 0$. Since $P$ is symplectically uniruled,  through  every  point $p\in P$ there exists  a $\sigma$-rational $J_P$ pseudo-holomorphic map $u$ passing through $p$. Since $J_P$ is fibered, the map $u_B:=\pi\circ u$ is $J_B$ pseudo-holomorphic.  Moreover, $u_B$ passes through $\pi(p)$. Consequently,  $B$ is covered by rational curves. Since $B$ is projective this implies that $B$ is symplectically uniruled; thus we have a contradiction.
\qed\\

Before proceeding to the applications, we make a remark concerning possible generalizations of Theorem \ref{mainobservation},  when we replace symplectic uniruledness by \emph{k-symplectic rational connectedness} with $k>1$ (cf Definition \ref{sympunirul}). Roughly speaking, a symplectic manifold is said to be $k$-symplectically rationally connected if there exists a non-vanishing Gromov-Witten invariant with $k$ points as constraints.  

   It is easy to find examples where Theorem \ref{SympDiv} and  observation \textbf{(O)} do not generalize to $k>1$, for instance when $P$ is a  trivial symplectic fibration.  Nevertheless,  it is possible to combine symplectic rational connectedness of the fiber and of the base in order to obtain symplectic rational connectedness of the total space.  Before making a precise  statement,  recall that the  relative  semi-positivity  condition on $(F,\om)$ implies that the fiber is  actually semi-positive. It is well-known that under such hypothesis, the \emph{Quantum homology} of $(F,\om)$ is generically well-defined. This quantum homology is endowed with a ring structure called \emph{quantum product}, which can be seen as a deformation of the intersection product in  homology. We introduce the following notation: let $[pt]$ be the   homology class of a point in $F$, we write $[pt]_Q$ to denote the class of the point seen as an element of the Quantum homology of $F$. Also, let $[pt]^k_Q$ denote the  $k$'th power of $[pt]_Q$ with respect to the quantum product.

\begin{theorem}\label{generalization} Suppose $(F,\om)$ is semi-positive relatively to $P$, and assume that $[pt]_Q^{l}\neq0$,  $l>1$. Furthermore, suppose that $(B,\om_B)$ is $(l+1)$-symplectically rationally connected for some class $\sigma_B\in H_2(B,\Z)$ admitting only simple decompositions. Then $P$ is $l$-symplectically rationally connected.
\end{theorem}

Next,  we give some applications of Theorem \ref{mainobservation} to the Weinstein Conjecture for separating hypersurfaces, or in shorter terms the \emph{$shW$-Conjecture}. 
 A. Weinstein conjectured in 1979 that:  "\emph{Every closed hypersurface $S$ of contact type in a given symplectic manifold $(X,\om)$ carries a closed characteristic}." \cite{Wein}. Since Viterbo's  proof of the conjecture for $(\R^{2n},\om_0)$ in 1986 \cite{V}, many results followed.  In particular, H. Hofer and C. Viterbo highlighted in \cite{HV} the strong interplay between genus zero Gromov-Witten invariants and the conjecture, in cases where the hypersurface $S$ \emph{separates} $X$, i.e. when there exist submanifolds $X^+$ and $X^-$ of $X$ having common boundary and such that $X=X^+\cup X^-$ and $S=X^+\cap X^-$. Note that this latter condition is realized whenever $H^1(X,\Z_2)=0$. Shortly after, these results were  extended by  G. Liu and G. Tian in \cite{LiuT}, to any symplectic manifold and any genus. In 2000, using the results of Liu and Tian, G. Lu proved the following:




\begin{theorem}\label{LuTh}(G. Lu \cite{Lu}, Corollary 3) Any separating hypersurface of contact type in a symplectically uniruled  closed symplectic manifold $(B,\om_B)$ admits a closed characteristic. 
In particular, the $shW$-Conjecture  holds in products  $(B\times F,\om_B\oplus\om)$ of  a symplectically uniruled  closed symplectic manifold $(B,\om_B)$ with any  symplectic manifold $(F,\om)$.
\end{theorem}

As a direct consequence of the first assertion of  Theorem \ref{LuTh} and of  Theorem \ref{mainobservation} we obtain the following generalization to  non-trivial Hamiltonian fibrations of the second assertion of Theorem \ref{LuTh}: 


\begin{cor}\label{WeinFib}  The  $shW$-Conjecture holds in $c$-splitting Hamiltonian fibrations with relatively semi-positive fiber and with  symplectically uniruled base for some class   admitting only simple decompositions. 
\end{cor}

As a consequence, we obtain the following which proof is given in Section 2: 
\begin{cor}\label{WeinFib2}  The  $shW$-Conjecture holds in Hamiltonian fibrations with relatively semi-positive fiber and with 2-symplectically rationnally connected base for some class   admitting only simple decompositions.
\end{cor}

In particular,  the $shW$-Conjecture holds in Hamiltonian fibrations (with relatively semi-positive fiber) over $(\C P^n,\om_{FS})$ where $\om_{FS}$ denotes the Fubini-Study Kahler form, and  over $(S^2\times S^2,\om\oplus\om)$ where $\om$ is an area form on $S^2$. It also holds for $c$-splitting Hamiltonian fibrations over ruled surfaces. Finally it should also hold for $c$-splitting Hamiltonian fibrations over blow-ups of such bases, since symplectic uniruledness is a symplectic cobordism invariant \cite{HLR}.\\

This note is organized as follows: in Section 1 we introduce the notions and notations needed for the proofs of the results, and in Section 2 we proceed with the proofs. \\
  
 \noindent\textbf{Aknowledgement.} I would like to thank Fran\c{c}ois Lalonde and Octav Cornea for their useful comments and  their encouragements. I would also like to thank Fran\c{c}ois   for his constant support and generosity. Finally, I thank R\'emi Lelclercq for his suggestions and comments on an earlier version of this note, that  helped improve the presentation of the paper.



\section{framework and tools}


\subsection{Gromov-Witten invariants} Let   $J$ be an $\om$-tame almost complex structure on $(X^{2n_X},\om)$, i.e. $J$ is a smooth endomorphism of $TX$ such that
 $$J^2=-Id_{TX} \hsp \text{and}\hsp\forall v\in TX, \om(v ,J v)>0.$$  
The set of all such almost complex structures is non-empty and contractible.  We  denote by $c_1^{TX}$ the associated first Chern class of $TX$, and say that $(X,\om)$ is \emph{semi-positive} if and only if there are no spherical class $A\in H_2(X;\Z)$ with positive  area $\om(A)$  such that 
\begin{equation*}3-n_X\leq c_1^{TX}(A)< 0.
\end{equation*}
We begin by recalling  the definition of genus zero Gromov-Witten invariant  for $(X,\om)$ assuming it is semi-positive. We follow  the expositions of  \cite{MS}, Chapter 7,  or \cite{RT}.  At the end we introduce the definition of simple decomposability.

\subsubsection{Gromov-Witten invariants} Recall that a  genus zero $J$-holomorphic map of $X$ is a smooth map $u:S^2\sra X$ solution to the Cauchy-Riemann equation:
 $$\delbar_Ju:=1/2(du+J\circ du\circ j_0)=0,$$
 where $j_0$ stands for the standard complex structure on $S^2\cong\C\cup\{\infty\}$.
Roughly speaking, a genus zero  Gromov-Witten invariant of $X$ is  a count, up to reparametrizations of the domain,  of  $J$-holomorphic rational maps  with marked points with values in $X$, satisfying some prescribed constraints at the marked points.  More precisely, consider the moduli space  of unparametrized genus zero $J$- holomorphic maps in $X$ with $l$ marked points and representing the class $A$
\begin{eqnarray*}\M_{0,l}(X,A,J)&:=& \Big\{(u,\mathbf{x}):=(u,x_1,...,x_l)\in C^{\infty}(S^2,X)\times (S^2)^l | x_i\neq x_j\hsp\text{if}\hsp i\neq j, \\
& &\,\,\,\,\left. \delbar_Ju=0\hsp,\hsp [u(S^2)]=A \Big\}\right/\G. 
\end{eqnarray*}
Let $\ov{\M}_{0,l}(X,A,J)$ denote its compactification, in the sense of Gromov, which coincides  with the moduli space of isomorphism classes $[\Sigma,u,\mathbf{x}]$ of stable $J$-holomorphic maps with $l$ marked points that represent $A$. 
This is a stratified space with finitely many strata. Its top stratum is the subset  $\M_{0,l}^{*}(X,A,J)\subset \M_{0,l}(X,A,J)$  consisting of  \emph{simple $J$-holomorphic maps}, i.e. maps that are somewhere injective. This stratum is, for generic $J$, a naturally oriented manifold with dimension  $2n_X+2c_1^{TX}(A)-6+2l$.  

There are two natural maps defined on $\ov{\M}_{0,l}(X,A,J)$, namely the \emph{evaluation at the marked points map}
$$ev^X_{l,J}:\ov{\M}_{0,l}(X,A,J)\sra X^l,\hspace{0.2cm} [\Sigma,u,x_1,\ldots,x_l]\mapsto (u(x_1),\ldots,u(x_l)),$$
and, for $l\geq3$, the \emph{forgetful-map map}
$$\mathfrak{f}:\ov{\M}_{0,l}(X,A,J)\sra \ov{\M}_{0,l}$$
assigning to every stable map the underlying reduced stable curve.
Formally,  when $l\geq3$, Gromov-Witten invariants are  the values of a multilinear homomorphism
\begin{equation*}\la\cdot\ra^X_{A,l}:H_*(\ov{\M}_{0,l},\Q)\otimes (H_*(X,\Q))^{\tens^{l}}\sra \Q
\end{equation*}
where
\begin{equation*}\la D;\beta_1,...,\beta_{l}\ra^X_{A,l}:=\int_{\ov{\M}_{0,l}(X,A,J)}(ev^X_{l,J})^*\left(PD_{X^l}(\beta_1\otimes ...\otimes \beta_l)\right)\cup \mathfrak{f}^*PD_{\ov{\M}_{0,l}}(D)
\end{equation*}
which is set to be zero unless:
\begin{equation*}\label{mathdim} 2n_X+2c_1^{TX}(A)=\sum_{i=1}^l (2n_X-\deg(\beta_i))-\deg(D).
\end{equation*} 
 When $l<3$, there is no forgetful map map and we simply "integrate" the   pull-back under $ev^X_{l,J}$ of the product of the $PD(\beta_i)$; in that case we use the notation $\la \beta_1,...,\beta_{l}\ra^X_{A,l}$, which could be viewed as  $\la [\ov{\M}_{0,l}];\beta_1,...,\beta_{l}\ra^X_{A,l}$ if we consider  $\ov{\M}_{0,l}$ as a manifold  of negative dimension $2l-6$. 

 The integration  above   has to be understood as evaluation of the cohomology class with the fundamental cycle of $\ov{\M}_{0,l}(X,A,J)$.  However, this space may not carry a fundamental class. In fact, lower strata in $\ov{\M}_{0,l}(X,A,J)$ may have  dimensions greater than   $\M_{0,l}^{*}(X,A,J)$, due to  the possible presence of stable maps with multiply covered components. This problem does not show up when the manifold is semi-positive. Concretely,  this condition imposes that the "boundary component" $\ov{\M}_{0,l}(X,A,J)\bs\M^*_{0,l}(X,A,J)$ is generically of   codimension at least two  with respect to  the top stratum. Then, consider cycles $V_1,...,V_l$ in $X$ representing the $\beta_i$'s, and a cycle $D$ in $\ov{\M}_{0,l}$  representing $D$.  Assume  that these cycles are in general position, and such that   
$\mathfrak{f}\times ev^X_{l,J}$ is strongly transverse to the product  $D\times V_1\times ...\times V_l$.  The corresponding \emph{Gromov-Witten invariant} is realized as the intersection number   $(\mathfrak{f}\times ev^X_{l,J}). (D\times V_1\times ...\times V_l)$ (which can be seen to be $\Z$-valued in that case). In particular, $GW$-invariants generically count simple maps.

\subsubsection{On the semi-positivity assumption.} We should mention that the semi-positivity assumption can be removed  using virtual perturbations of the Cauchy-Riemann equation as in \cite{LT},\cite{RT},\cite{Rvirt}, amongst others. However, we will not work in such generality. Also, let us emphasize that the semi-positivity assumption can   be dropped if   the boundary of $\ov{\M}_{0,l}(X,A,J)$ only consists of \emph{simple $J$-holomorphic stable maps}. Namely, we say that a stable $J$-holomorphic map is \emph{simple} if and only if it has  no non-constant multiply covered component and no two non-constant components having the same image in $X$.  This leads to the following definition,

\begin{defn}\label{irred} Let $A\in H_2(X;\Z)$ be a spherical  class representable by $J$-holomophic maps.  We say that $A$  \emph{only admits simple decompositions} if  every stable $J$-holomorphic  map  in  $\ov{\M}_{0,l}(X,A,J)$ has no non-constant multiply covered component and no two non-constant components having the same image in $X$. 
\end{defn}

This is what we request from  $\sigma_B$ in Theorem \ref{mainobservation}.   Note that the set of $\om$-tame almost complex structures $J$ with respect to which a given class $A$
only admits simple decompositions is open in the set of $\om$-tame almost complex structures. However,  we cannot make sure that this set  is connected or non-empty.

\subsection{Symplectic rational connectedness.} Using the notations introduced  in the preceding paragraph,  we   define  the notion of  \emph{$k$-symplectic rational connectedness}, following \cite{HLR} and \cite{Lu}. 

\begin{defn}\label{sympunirul} Let $k>0$ be an integer, and let $\sigma\in H_2(P;\Z)$ be a spherical homology class. A symplectic manifold $(X,\om)$ is \emph{$k$-symplectically rationally connected for $\sigma$}, or simply \emph{k-SRC} for $\sigma$,  if and only if there exist classes   $\beta_{k+1},...,\beta_l\in H_*(X,\Q)$, and $D\in H_*(\ov{\M}_{0,l},\Z)$ such that:
\begin{equation*}\la D; [pt],...,[pt],\beta_{k+1},...,\beta_l\ra_{A,l}^X\neq 0.\end{equation*}
we will say that $(X,\om)$ is \emph{$k$-symplectically rationally connected} if there exists $\sigma$ such that it is \emph{k-SRC} for $\sigma$. 
If  $k=1$, we will say that  $(X,\om)$ is \emph{symplectically uniruled} or \emph{SU}.
Furthermore, if $(X,\om)$ is symplectically uniruled and $l=3$, we say that $(X,\om)$ is \emph{strongly symplectically uniruled} or \emph{SSU}. 
\end{defn}

In particular, $(X,\om)$ is $k$-SRC for $A$ only if, through every $k$ generic points of $X$, there is a genus zero pseudo-holomorphic map representing $A$. The converse may not be true in general. Nevertheless, the equivalence holds  in smooth projective varieties. 

\begin{theorem}(Ruan \cite{Rvirt},Koll\`{a}r \cite{Kollar1})\label{unirulingproj} A smooth projective variety is symplectically uniruled if and only if it is \emph{projectively uniruled}, i.e. through every point of the manifold there is a holomorphic map.  
\end{theorem}

  Actually, it follows from the \emph{splitting axiom for Gromov-Witten invariants} that  a projective manifold is (projectively) uniruled if and only if it is strongly symplectically uniruled, as pointed out by Ruan \cite{Rvirt}. The main ingredient of the proof is the following property of rational curves in projectively uniruled manifolds $(X,J,\om)$  due to J. Koll\`{a}r, Y. Miyaoka, and S. Mori (see Koll\`{a}r \cite{Kollar}, Theorem 3.11):
  for a very general point $b\in X$, if $u:\C P^1\sra X$ is  a morphism such that $[u(\C P^1)]\neq 0$ and $u(0)=b$, then $H^1(\C P^1,u^*TX)=0$.   Such morphism is said to be \emph{free (over 0)}.  Now, for $b\in X$ and a spherical class $\sigma$, let $\M(X,\sigma,J;b)$ denote the moduli space of $\sigma$ rational $J$-holomorphic maps  $u:\C P^1\sra X$ such that $b\in \text{Im}(u)$. We will say that $\sigma$ is \emph{free} if, for every general point $b\in X$, the moduli space  
  $\M(X,\sigma,J;b)$  only consists of free morphisms.   For our purpose, we  mildly refine  the statement of Theorem \ref{unirulingproj} by the following straightforward observation.

  \begin{lem}\label{projuni} If  a smooth projective variety $(X,J,\om)$ is   uniruled, it is strongly  symplectically uniruled for a  free class $\sigma$ admitting only simple decompositions.
  \end{lem}

\proof Since $X$ is uniruled, for every $b\in X$ there exists a morphism $u:\C P^1\sra X$ such that $b\in \text{Im}(u)$. Now, fix a sufficiently general point $b$. In the proof of Theorem \ref{unirulingproj} given in \cite{Rvirt},  Ruan shows the equality between
\begin{eqnarray*}N_1:=\min\left\{\om(\sigma)>0| \M(X,J,\sigma;b)\neq \emptyset\right\} 
\end{eqnarray*}
and 
\begin{eqnarray*}N_2:=\min\left\{\om(\sigma)>0| \exists \alpha_1,...,\alpha_l\in H_*(X), \la [pt];[pt],\alpha_1,...,\alpha_l\ra^X_{\sigma,l+1}\neq 0\right\}.
\end{eqnarray*}
He also points out that  $N_2$ is realized by classes with non trivial three point $GW$-invariant. Suppose $\sigma\in H_2(X;\Z)$ realizes this minimum. Then every rational curve passing through $b$ and representing $\sigma$ is irreducible, i.e. any stable holomorphic map through $b$ is simple. It follows that   $\M(X,J,\sigma;b)$ is   compact. Furthermore,   by choosing $b$ general enough,  for any $u\in \M(X,J,\sigma;b)$ the obstruction $H^1(\C P^1,u^*TX)$ vanishes; thus, the moduli space $\M(X,J,\sigma;b)$ is a smooth oriented manifold, as desired.
\qed


 \subsection{Hamiltonian fibrations and the product formula} In  this section we  introduce the notion of Hamiltonian fibrations which provides natural framework to study Gromov-Witten invariants. 
 
\subsubsection{Hamiltonian fibrations.}   By definition, a \emph{symplectic  fibration} is a locally trivial smooth fibration with symplectic reference fiber $(F,\omega)$,
$$(F,\om)\stackrel{\iota}{\hra} P\stackrel{\pi}{\sra} B,$$  
 and which structure  group lies in the group of symplectic  diffeomorphisms of the fiber, denoted $\mathrm{Symp}(F,\omega)$. It follows that each  fiber  $F_b:=\pi^{-1}(b)$  is naturally  equipped with a symplectic form $\omega_b$. A symplectic fibration is  \emph{Hamiltonian} if the structure group can be reduced to the group $\mathrm{Ham}(F,\omega)$ of Hamiltonian diffeomorphisms. Hamiltonian fibrations are characterized as follows:
 \begin{theorem}(\cite{MS2}, Theorem 6.36.) A symplectic fibration $P$ as above is Hamiltonian if and only if the following two conditions are verified:
 \begin{itemize}
\item[($H_1$)]  $P$ is symplectically trivial over the $1$-skeleton of $B$;
\item[($H_2$)]  there exists a unique closed  connection 2-form $\tau\in \Omega^2(P)$ extending   the family $\{\omega_b\}_{b\in B}$  such that  the integration of $\tau^{n_F+1}$ over the fibers of $P$ vanishes.
\end{itemize}
\end{theorem}
 
 The closed $2$-form in the theorem above is usually refered to as the \emph{coupling form}.  The coupling form  defines a connection on $P$, i.e. a splitting at each $p\in P$,
 $$T_pP=Hor_{\tau,p}\oplus \ker d\pi(p).$$
The theorem  states  that the corresponding holonomy  around any  loop in $B$ is in $\mathrm{Ham}(F,\omega)$. Any other closed extension  $\tau'$ of $\om$ generating the  same  horizontal distribution $Hor_{\tau}$   is actually uniquely obtained from  $\tau$  via the equation:
$$\tau'=\tau+\pi^*\varrho,\hsp\varrho\in\Omega^2(B).$$   
One of the main features of Hamiltonian fibrations over closed symplectic bases  is that  they can be given symplectic structures compatibly with the family  $\{\omega_b\}_{b\in B}$.
Such symplectic structures are  obtained as follows: 
$$\omega_{P,\kappa}:=\tau+\kappa \pi^*\omega_B,$$ where $\kappa>0$ is a large enough real number  such that $\omega_{P,\kappa}$ is non-degenerate. Hence, this class of fibrations provide a nice framework to define Gromov-Witten invariants.
We will also use   $\tau$   to denote the   deRham cohomology class corresponding to $\tau$. There is another canonical cohomology class of $P$ that will play an important part in the next paragraphs, namely the vertical Chern class $c_v\in H^2(P;\Z)$. This class is defined  as the first Chern class of the vertical subbundle $\ker d\pi$. 



 
As a specific example, let us mention  the case of Hamiltonian fibrations over $S^2$ with fiber $(F,\om)$, which will play an important part in this note.  This  class of   examples is particularly important in  symplectic topology  due to the correspondence between these fibrations and the fundamental group of  $Ham(F,\om)$. More precisely, for $\gamma\in\pi_1(Ham(F,\om)) $ one defines a Hamiltonian bundle $\pi:P_{\gamma}\sra S^2$ with fiber $F$ via the clutching construction: choose any representative $\tilde{\gamma}:[0,1]\sra Ham(F,\om)$ for $\gamma$, then
$$P_{\gamma}:=\frac{(D^+\times F)\sqcup (D^-\times F)}{(e^{2\pi i\theta},x)\sim (e^{-2\pi i\theta},\wt{\gamma}(\theta).x),\hsp\text{on}\hsp S^1\times F}$$
 where $D^{\pm}\subset \C$ denotes the closed  unit disc. This construction is independent of the representative $\wt{\gamma}$;  moreover, any Hamiltonian fibration over $S^2$ with $(F,\om)$ as fiber can be constructed in this way   (see \cite{Se} or \cite{LMP}). In this context, the coupling class  will be denoted by $\tau_{\gamma}$  and  $c_{\gamma}$ will denote the  vertical Chern class associated to $P_{\gamma}$.

 \subsubsection{The product formula.} Following \cite{H}, we give a product formula for $GW$-invariants of a Hamiltonian fibration, assuming the fiber  $(F,\om)$ is  \emph{semi-positive relative to the total space}, i.e.   
 \begin{equation}\tag{$\star$}\label{ssp}\forall A\in H_2^S(F):\hsp  \omega(A)>0,\hsp c^v(\iota(A))\geq 3-n_P\hspace{0.3cm}\Longrightarrow c^v(\iota(A))\geq 0, 
\end{equation}
 where $H_2^S(F)$ is the spherical homology subgroup of $H_2(F,\Z)$ (i.e.  the image of $\pi_2(F)$ under the Hurewicz map) and $\iota$ denotes the map in homology induced from the natural embedding of the fiber. For instance, this implies  that the fiber is semi-positive.
For this purpose we equip  $P$ with an \emph{(almost) complex structure} $J_P$ which is compatible with the fibration structure and a Hamiltonian connection in the sense given below:
\begin{defn}\label{fiberedacs} An almost complex structure $J_P$ on $P$ is said to be  \emph{compatible with $\pi$ and $\tau$}, or just \emph{fibered}, if and only if there exists   an $\om_B$-tame complex structure $J_B$ on $B$ and a family of  $\om_b$-tame almost complex structures $J_b$ in $F_b$ such that:
\begin{itemize}
  \item $d\pi \circ J_P=J_B\circ d\pi$,
  \item $J_b:=\left.J_P\right|_{F_b}$ for all $b\in B$,
  \item $J_P$ preserves the horizontal distribution induced by $\tau$.
\end{itemize}
\end{defn} 
For fibered $J_P$,  the projection $\pi$ induces a map between moduli spaces:
$$ \ov{\pi}:\MMC_{0,l}(P,\sigma,J_P)\sra\MMC_{0,l}(B,\sigma_B,J_B),\hsp [\Sigma,u,\mathbf{x}]\mapsto [\Sigma, \pi(u),\mathbf{x}]$$
where $\sigma_B:=\pi_*\sigma$. In what follows we  assume that $\sigma_B$  is non-zero.  The fiber of $\ov{\pi}$ over $[\Sigma, u_B,\mathbf{x}]$ can be described as follows.  
Let $C$ denote the image of $u_B$ in $B$, and let $P_C$ denote the restriction of $P$ along $C$; $P_C$ defines a Hamiltonian fibration  over $S^2$ with   coupling form given by the  pull-back of  $\tau$ under the natural inclusion $\iota_{P_C}:P_C\hra P$.  If  $J_C$ denotes the fibered almost complex structure on $P_C$ given by the restriction of $J_P$ to $P_C$ we have the following identification:
$$\ov{\pi}^{-1}([\Sigma,u_B,\mathbf{x}])\equiv\bigsqcup_{B_{\sigma}:=\{\sigma'\in H_2(P_C;\Z)|\iota_{P_C}\sigma'=\sigma\}} \left(\M_{0,l}(P_C,J_C, \sigma')\cap \mathfrak{f}^{-1}([\Sigma,\mathbf{x}])\right).$$
Regarding evaluation at the marked points maps we  have the commutative diagram:
\begin{equation}\label{evdiagram}
\xymatrix{\ov{\pi}^{-1}([\Sigma,u_B,\mathbf{x}]) \ar[r]\ar[d]^{ev_{(u_B,\mathbf{x})}}&\MMC_{0,l}(P,\sigma,J_P)\ar[d]^{ev^P_{l,J_P}}\ar[r]^{\ov{\pi}}& \MMC_{0,l}(B,\sigma_B,J_B) \ar[d]^{ev^B_{l,J_B}}\\
F^l\ar[r]^{(\iota)^l}&P^l\ar[r]^{\pi^l} & B^l}
\end{equation}
where 
$$ev_{(u_B,\mathbf{x})}:\ov{\pi}^{-1}(u_B,\mathbf{x})\sra F^l,\hspace{0.5cm} u\mapsto (u(x_1),...,u(x_l))\in\prod_{i=1}^lF_{u_B(x_i)}.$$
The product formula is obtained by considering the (respective) intersections of $ev_{(u_B,\mathbf{x})}$, $ev^B_{l,J_B}$, and $ev^P_{l,J_P}$, with the  product cycles: 
\begin{equation*} \mathcal{C}^F:=\prod_{i=1}^{l}V^F_i,\hspa\mathcal{C}^B:=\prod_{i=1}^{l}V^B_i,\hspa \mathcal{C}^P:=\prod_{i=1}^{l}V^P_i,
\end{equation*}
where, $V^F_i$, $V^B_i$, and $V^P_i$,  respectively represent   homology classes,  $c_i^F$, $c_i^B$, and $c_i^P$,  verifying that for some integer $0\leq m\leq l$:
\begin{equation}\tag{$\star\star$}\label{cond11}
\begin{cases} c_i^B=pt, \hsp c_i^P=\iota(c_i^F)  & \text{for $i=1,...,m$} \\
 c_i^F=[F],\hsp c_i^P=\pi^!(c_i^B)& \text{for $i=m+1,...,l$.}
\end{cases}
\end{equation}
 where $\pi^!$ stands for the shriek map:
\begin{equation*}\pi^{!}:H_*(B)\sra H_{2n_F+*}(P), \hsp \alpha\mapsto PD_P^{-1}\pi^*PD_B(\alpha).
\end{equation*}
In other words, we consider cycles  in $P$ that are either  images of cycles in $F$ under $\iota$, or preimages under $\pi$ of cycles in $B$. Let $\iota_C$ denote   the map in homology induced from the natural inclusion of $F$ in $P_C$. We have,

 \begin{theorem}(\cite{H}, Theorem A)\label{productformula0} 
Let $\pi: P\sra B$ be a Hamiltonian fibration  with relative semi-positive fiber $(F,\om)$. Let $\sigma\in H_2(P,\Z)$ and 
suppose  $\sigma_B:=\pi_*(\sigma)\neq 0$ only admits simple decompositions for some $\om_B$-tame almost complex structure $J_B$ in $B$. 
Let  $c_i^P,c_i^B,c_i^F$ be as in \eqref{cond11}. For generic fibered almost complex structure  lifting $J_B$, we  have:
\begin{equation}\label{simpversionPF}\la D;c^P_1,...,c_l^P\ra_{\sigma,l}^P=\la D;c^B_1,...,c_l^B\ra_{\sigma_B,l}^B\cdot\sum_{\sigma'\in B_{\sigma}} \la [pt]; \iota_C(c^F_1),...,\iota_C(c^F_l)\ra_{\sigma',l}^{P_C}
\end{equation}
where $C$ is a curve counted in $\la D; c_1^B,...,c_l^B\ra^B_{\sigma_B,l}$ and $D\in H_*(\ov{\M}_{0,l})$.
 \end{theorem}

This formula in particular states that the sum in the right handside of the formula  is independent of  the chosen $C$. Note that this sum is well-defined due to Gromov 's compactness.  

  The core of the proof consists in  establishing that the  $GW$-invariants  involved are generically and simultanuously well-defined, in other words the problem is to realize transversality while preserving   the map $\ov{\pi}$ defined above. The proof of this is based on the relation:
$$\pi_*\circ D^P=D^B\circ \pi_*,$$
where $D^P$ and $D^B$ respctively stand for the Fredholm operators obtained by linearizing the Cauchy Riemann operators  $\delbar_{J_P}$ and $\delbar_{J_B}$. As a consequence, we derive  an  exact sequence
\begin{equation*}
  0\rightarrow \ker D^{v} \rightarrow \ker D^P \rightarrow \ker D^B\rightarrow \coker D^{v}\rightarrow \coker D^P\rightarrow\coker D^B\rightarrow 0,
\end{equation*}
where  $D^v$ denotes  the restriction of $D^P$ to vector fields along the curves that are vertically valued, i.e. with values in $\ker d\pi$. The vanishing, at least at the level of the universal moduli spaces, of the obstructions in the sequence above is provided by: 1)  the irreducibility hypothesis on $\sigma_B$ for the vanishing of the last term of the sequence; 2) perturbing the Hamiltonian connection for the vanishing of $\coker D^v$.  It follows from standard arguments that  for generic fibered almost complex structure the following holds:
\begin{itemize}
\item the subset $\MMC^{**}_{0,l}(P,\sigma,J_P)$ of $\MMC_{0,l}(P,\sigma,J_P)$ consisting of simple maps that project to simple maps under $\ov{\pi}$ and the moduli space $\MMC^*_{0,l}(B,\sigma_B,J_B)$ 
are oriented manifolds; 
\item
for countably many  $(u,\mathbf{x})\in \MMC^*_{0,l}(B,\sigma_B,J_B)$,  the preimage $\ov{\pi}^{-1}(u,\mathbf{x})$ is an   oriented manifold.
\end{itemize}

More generally, $\ov{\pi}$ extends to a map, using stabilization, between the compactifications $\ov{\MMC}_{0,l}(P,\sigma,J_P)$ and $\ov{\MMC}_{0,l}(B,\sigma_B,J_B)$. We can repeat the arguments above for each stratum in $\ov{\MMC}_{0,l}(P,\sigma,J_P)$ projecting to some stratum in  $\ov{\MMC}_{0,l}(B,\sigma_B,J_B)$. The   irreducibility hypothesis on the decompositions of  $\sigma_B$ ensures that the $\coker D^B$ term always vanishes and  that the image under the evaluation map $ev^B_{l,J_B}$ of the lower strata in  $\ov{\MMC}_{0,l}(B,\sigma_B,J_B)$  have  codimension at least 2 with respect to the top stratum $\M^*(B,\sigma_B,J_B)$. Furthemore, condition \eqref{ssp} ensures  that the image under the evaluation map $ev^P_{l,J_P}$ of the lower strata in  $\ov{\MMC}_{0,l}(P,\sigma,J_P)$  have  codimension at least two with respect to $\MMC^{**}_{0,l}(P,\sigma,J_P)$. Once transversality is achieved, one recovers the formula using  diagram \eqref{evdiagram} and the following observations:
\begin{itemize}
\item  since all cokernels vanish, it follows from the exact sequence above that:
$$\det (\ker D^P)\cong \det (\ker D^B) \otimes \det (\ker D^v).$$
Thus,  the orientation of  $\MMC^{**}_{0,l}(P,\sigma,J_P)$ is  given by the product of the orientations of $\M^*_{0,l}(B,\sigma_B,J_B)$ and   $\ov{\pi}^{-1}(u_B,\mathbf{x})$,  where  $(u_B,\mathbf{x})\in \M^*_{0,l}(B,\sigma_B,J_B)$. Moreover, the orientations of the product pseudo-cycles are also given by a product:
\begin{equation*}
 \det T\CC_P\cong \det T\CC_B\otimes \det T\CC_F.
\end{equation*}
 \item It follows from symplectic triviality of $P$   over  the $1$-squeleton of $B$, that for any two $J_B$-holomorphic maps representing the same  class $\sigma_B$, the restrictions of $P$ to the images of the two maps are isomorphic as Hamiltonian fibrations. Hence, $\la\phantom\cdots\ra^{P_C}_{\sigma',l}$ does not depend on $C$.
\end{itemize}
 The  proof is concluded by showing, using standard arguments,  the independence with respect to the generic fibered almost complex structure of the Gromov-Witten invariants involved in the formula. See \cite{H} for the details.

\subsection{Quantum Homology and Seidel elements} We start with the definition of (small) Quantum homology with universal Novikov ring. As  a module the   Quantum homology of $(X,\om)$ is given by:
\begineqno QH_*(X)\equiv QH_*(X,\Lambda):=H_*(X,\Q)\otimes \Lambda
\endeqno
where $\Lambda$ denotes  some coefficient ring that we specify.  We will take
$\Lambda:=\Lambda^{univ}[q^{-1},q]$
where $\Lambda^{univ}$ is the ring of generalized Laurent series in variable $t^{-1}$, i.e. an element $\lambda\in\Lambda^{univ} $ can be written as a formal sum
\begineqno \lambda=\sum_{i\geq 0} \lambda_it^{r_i},\hspa \lambda_i\in \Q,r_i\in \R,\hsp r_i>r_{i+1},\hsp \lim_{i\sra\infty}r_i=-\infty.
\endeqno
The grading on $\Lambda$ is given by imposing that $\deg(q)=2$. Let $\Lambda_j$ denote the set of elements of degree $2j$ in $\Lambda$; then we give  quantum homology the following grading:
\begineqno QH_k(X):=\bigoplus_{i+2j=k} H_i(X,\Q)\otimes \Lambda_j
\endeqno 
Next, we introduce the (small) quantum product:
$$\star :QH_i(X)\otimes QH_j(X)\sra QH_{i+j-2n}(X),\hsp (a,b)\mapsto a\star b.$$ 
First, let $\{e_{\nu}\}$ be a basis of homology of $H_*(X)$, and let $\{e^*_{\nu}\}$ denote the corresponding dual basis with respect to the intersection pairing.  
Then, the quantum product is defined as follows: for $a\in H_i(X,\Q)$ and $b\in H_j(X,\Q)$, set
\begineq a\star b=\sum_{B\in H^S_2(X),\nu}\la a,b,e_{\nu}\ra^X_{B,3} e^*_{\nu}\otimes q^{-c_1(B)}t^{-\om(B)}
\endeq
We extend this by linearity with respect to $\Lambda$. Note that $[X]=1$ is the identity for this ring structure. Now, consider the  $\Lambda$-submodule:
\begineqno \mathcal{Q}_-:=\bigoplus_{i<2n} H_i(X)\otimes \Lambda.
\endeqno
The following lemma of McDuff (\cite{MDuniruled}, Lemma 2.1) will be  essential  when proving Theorem \ref{mainobservation}. We give its proof for the reader's convenience:

\begin{lem}\label{explicitSeidelelement} If $(X,\om)$ is not strongly uniruled then $\mathcal{Q}_-$ is an ideal in $QH_*(X)$. Furthermore, if $\mathcal{Q}_-$ is an ideal  $a\in QH_*(X)$ is a unit, then there exists $x\in\mathcal{Q}_-$ and $\lambda\neq 0$ in $\Lambda$ such that: 
$$a=1\otimes\lambda+x.$$
\end{lem}
 
\proof Let $c\in  QH_*(X)$ and assume by contradiction that there exists $b\in\mathcal{Q}_-$ such that $c\star b\notin \mathcal{Q}_- $. By definition of the quantum product this means that there exists $B\in H_2^S(X)$ such that 
$$\la c,b,pt\ra^X_{B,3}\neq 0;$$
hence $(X,\om)$ is strongly uniruled which is impossible. Next, if $a$ is a unit, then $a\notin \mathcal{Q}_-$, for if it was we would have that $1$ belongs to $\mathcal{Q}_-$. Therefore, we can write 
$a=1\otimes\lambda+x$ with $x\in\mathcal{Q}_-$ and $\lambda\neq 0$.
\qed\\

\noindent\emph{Seidel's representation.}  This is a representation of  $\pi_1(Ham(X,\om))$ in the subring  $QH^{\times}_*(X,\Lambda)$ of units of  the quantum homology of $(X,\om)$. Recall that $\pi:P_{\gamma}\sra S^2$ denotes the Hamiltonian fibration obtained from  $\gamma\in\pi_1(Ham(X,\om)) $  via the clutching construction. Also, let $\H_{\gamma}\subset H^S_2(P_{\gamma})$ denote the subset of section classes, i.e. of spherical classes that project to $[S^2]$ under the fibration projection.  

\begin{defn} The Seidel representation map
$$S:\pi_1(Ham(X,\om))\sra QH^{\times}_*(X,\Lambda),\hsp \gamma\mapsto S(\gamma)$$
is defined as follows:
\begin{equation} S(\gamma):=\sum_{\sigma\in \H_{\gamma} ,\nu}\la \iota(e_{\nu})\ra^{P_{\gamma}}_{\sigma,1} e^*_{\nu}\otimes q^{-c_{\gamma}(\sigma)}t^{-\tau_{\gamma}(\sigma)}
\end{equation}
\end{defn}
Geometrically,  $S(\gamma)$ "counts"    holomorphic sections of $P_{\gamma}$ intersecting the cycles $e_{\mu}$ in the  fiber above  the north pole of $S^2$. The \emph{splitting axiom}, \cite{MS} Theorem 11.4.1, for fibrations over $S^2$  gives that for all section classes $\sigma \in H_2^S(P_{\gamma})$ and for classes $\alpha_1,...,\alpha_l\in H_*(F)$ the following holds for every integer $0\leq k\leq l$:
\begin{eqnarray*}
& &\la[pt];\iota(\alpha_1),...,\iota(\alpha_l)\ra^{P_{\gamma}}_{\sigma,l} =
\\&& \sum_{A\in H_2^S(F;\Z),\nu}
\la[pt];\iota(\alpha_1),...,\iota(\alpha_{k}),\iota(e_{\nu})\ra^{P_{\gamma}}_{\sigma-\iota(A),k+1}  
\la[pt];e^*_{\nu},\alpha_{k+1},...,\alpha_l\ra^{F}_{A,l-k+1} 
\end{eqnarray*}
It follows easily  that  the $\Lambda$ linear action of $S(\gamma)$ on $a\in H_*(F)$ is given by:
\begin{equation}\label{Seidelformula} S(\gamma)(a):=S(\gamma)\star a=\sum_{\sigma\in \H_{\gamma} ,\nu}\la \iota(a), \iota(e_{\nu})\ra^{P_{\gamma}}_{\sigma,2}e^*_{\nu}\otimes q^{-c_{\gamma}(\sigma)}t^{-\tau_{\gamma}(\sigma)}
\end{equation}
Note that in the above notations $S(\gamma)=S(\gamma)(1)$. Let us define the following equivalence class on section classes: we say that $\sigma_1$ is equivalent to $\sigma_2$ if and only if $$\tau_{\gamma}(\sigma_1-\sigma_2)=0=c_{\gamma}(\sigma_1-\sigma_2).$$
We will use the notation $[\sigma]$ to denote the equivalence class of $\sigma$. From Theorem 3.A in \cite{LMP}, we have 

\begin{lem}\label{LemmeSeidelinverse}  For every $\gamma\in \pi_1(Ham(F,\om))$ and every non zero $a\in H_*(F)$, there exists an equivalence class of section classes $[\sigma]$ and $b\in H_*(F)$ such that:
 $$\sum_{\sigma'\in [\sigma]}\la\iota (a),\iota (b) \ra^{P_{\gamma}}_{2,\sigma'}\neq 0.$$
 In particular there is a section class $\sigma\in \H_{\gamma}$ such that $\la\iota(a),\iota(b) \ra^{P_{\gamma}}_{2,\sigma}\neq 0$. 
\end{lem}

\proof   Assume it is  not true. Equation \eqref{Seidelformula}  and linearity of Gromov-Witten invariants then imply that $S(\gamma)(a)=0$. But since $S(\gamma)$ is invertible, this is only possible for $a=0$; hence the contradiction.  \qed \\

\section{proofs of the results} 
We are now ready to prove the results. We begin by proving Theorem \ref{mainobservation} and Theorem \ref{generalization}.   More precisely,  we prove the following:

\begin{theorem}\label{unirulingfib2} Assume $(F,\om)$ verifies \eqref{ssp} and is not strongly uniruled. If $(B,\om_B)$ is symplectically uniruled for some class $\sigma_B$ admitting only simple decompositions, then $P$ is also symplectically uniruled. Moreover, if   $[pt]_Q^{l}\neq0$,  and if $(B,\om_B)$ is $(l+1)$-SRC for   $\sigma_B$ admitting only simple decompositions,  then $P$ is at least $l$-SRC.
\end{theorem}

\proof Since $B$ is symplectically uniruled there exists  $\sigma_B\in H_2^S(B)$ and classes $c_1^B,...,c^B_l\in H_*(B)$ such that:  \begineq\label{uniB} \la D; [pt],c_1^B,...,c^B_l\ra_{\sigma_B,l+1}^B\neq 0.
\endeq
 Let $C$ be the image of a  map  counted in \eqref{uniB}.  The restriction $P_C$ of $P$ to $C$ is a Hamiltonian fibration over $S^2$.  Let $\phi\in \pi_1(Ham(F,\om))$  be a Hamiltonian loop corresponding to this fibration, and let $S(\phi)\in QH^{\times}(F;\Lambda)$ denote the corresponding Seidel element. Since  $(F,\om)$ is not strongly uniruled,  it follows from  Lemma \ref{explicitSeidelelement} that there exists $0\neq \lambda\in \Lambda$ and $x\in \mathcal{Q}_- $ such that:
$$S(\phi)=1\otimes \lambda+x.$$
In particular, this directly implies that there is an equivalence class  $[\sigma]$ of section classes  in $P_C$ such that
$$0\neq \sum_{\sigma'\in[\sigma]} \la[pt]\ra^{P_C}_{\sigma',1}=\sum_{\sigma'\in[\sigma]}\la[pt]; [pt],\iota_C([F]),...,\iota_C([F])\ra^{P_C}_{\sigma',l+1}. $$

Let $\iota_{P_C}$ denote the inclusion in homology induced by the inclusion of $P_C$ into $P$. It is easy to see that the image under $\iota_{P_C}$ of $[\sigma]$ defines an equivalence class of spherical classes projecting on $\sigma_B$, in the sense that for any $\sigma_1,\sigma_2\in \iota_{P_C}([\sigma])$ we have $$\tau(\sigma_1-\sigma_2)=0=c_v(\sigma_1-\sigma_2).$$
 By the product formula for Gromov-Witten invariants \eqref{simpversionPF}:
\begin{eqnarray*}&&\sum_{\sigma'\in \iota_{P_C}([\sigma])}\la D; [pt],\pi^{!}(c_1^B),..., \pi^{!}(c_l^B)\ra_{\sigma',l+1}^P\\
&&=\sum_{\sigma'\in  \iota_{P_C}([\sigma])}\la D; [pt], c^B_1,...,c_l^B\ra_{\sigma_B,l+1}^B\sum_{\sigma''\in B_{ \sigma'}}\la[pt]; [pt],\iota_C([F]),...,\iota_C([F])\ra^{P_C}_{\sigma'',l+1}.\\
&&=\la D; [pt], c^B_1,...,c_l^B\ra_{\sigma_B,l+1}^B \sum_{\sigma'\in[\sigma]}\la[pt]; [pt],\iota_C([F]),...,\iota_C([F])\ra^{P_C}_{\sigma',l+1}\neq 0\\
\end{eqnarray*}
In particular,  there is at least one $\sigma'\in \iota_{P_C}([\sigma])$ such that $$\la D; [pt],\pi^{!}(c_1^B),..., \pi^{!}(c_l^B)\ra_{\sigma',l+1}^P\neq 0.$$

Now we prove the second assertion of the theorem.  Since for every loop $\phi$ of Hamiltonian diffeomorphisms of $F$ the Seidel element $S(\phi)$ is an unit in $QH_*(F;\Lambda)$, we have that
$$0\neq S(\phi)\star ([pt]_Q^l).$$
Using the splitting axiom for Gromov-Witten inviariants, and from the definition of the Seidel element,  one obtains that
\begin{eqnarray*} 
S(\phi)\star ([pt]_Q^l)&=&\sum_{ A\in H_2^S(F),
    \nu}\la [pt]; [pt],...,[pt], e_{\nu}\ra^{F}_{A,l+1}S(\phi)(e^*_{\nu})\otimes q^{-c^{TF}_1(A)}t^{-\om(A)}\\
&=& \sum_{\sigma\in \H_{\phi},
 A\in H_2^S(F),
    \nu,\mu}\la [pt]; [pt],...,[pt], e_{\nu}\ra^{F}_{A,l+1}  \la \iota(e^*_{\nu}),\iota(e_{\mu})\ra^{P_{\phi}}_{\sigma,2} \\
    &&\hspace{3.5cm}e^*_{\mu}\otimes q^{-c_{\phi}(\sigma+\iota(A))}t^{-\tau_{\phi}(\sigma+\iota(A))}  \\
& =&\sum_{
\wt{\sigma}\in \H_{\phi},
    \mu}\la [pt]; [pt],...,[pt], \iota(e_{\mu})\ra^{P_{\phi}}_{\wt{\sigma},l+1}e^*_{\mu}\otimes q^{-c_{\phi}(\wt{\sigma})}t^{-\tau_{\phi}(\wt{\sigma})}
\end{eqnarray*}
 Hence,  there is a class $a\in H_*(F)$ and a class of section classes $[\sigma]$ such that
$$\sum_{\sigma'\in[\sigma]}\la[pt]; [pt],...,[pt],\iota(a)\ra^{P_{\phi}}_{\sigma',l+1}\neq 0.$$
We conclude by the use of the product formula as before. 
\qed\\

Now Theorem  \ref{mainobservation} follows easily:\\

\noindent\textbf{Proof of Theorem \ref{mainobservation}:} In case where $(F,\om)$ is strongly uniruled, the conclusion follows from Theorem \ref{SympDiv}. Otherwise, we simply use Theorem \ref{unirulingfib2} to conclude.
\qed\\

 As for Theorem \ref{mainobservation}, Corollary  \ref{cormainobs} follows directly from Theorem \ref{SympDiv} and the proposition below. This result is a simple consequence of Theorem \ref{unirulingproj} asserting that projectively uniruled manifolds are strongly symplectically uniruled.  

\begin{prop}Assume $(F,\om)$ verifies \eqref{ssp} and is not strongly uniruled. Also assume that $(B,J_B,\om_B)$ is  a  uniruled projective manifold. Then $P$ is SSU. 
\end{prop}

\proof From lemma \ref{projuni}, $(B,J_B,\om_B)$ is uniruled for a class  $\sigma_B$ verifying the assumptions of Theorem \ref{unirulingfib2}, and which is obstruction free. Hence, the product formula can be applied and the proof follows.
\qed\\

Now, we proceed to the proof of Corollary  \ref{WeinFib2}:\\
 
\noindent\textbf{Proof of Corollary \ref{WeinFib2}:}  By definition, $B$ is 2-SRC. We will show that $(F,\om)$ homologically injects in $P$.  Assume this is not true, then there exists a non-zero element $\ov{a}\in H_*(F,\Z)\cap \ker\iota$. It follows by linearity that any Gromov-Witten invariant with entry  $\iota(\ov{a})$ must vanish. Now, let $C$ be the image of a  map  counted in $$\la D;[pt],[pt],c_3^B,...,c^B_l\ra_{\sigma_B,l}^B\neq 0.$$
Since $P_C$ is a Hamiltonian fibration over $S^2$, we can apply Lemma \ref{LemmeSeidelinverse}, taking  $a=[pt]$, and find that 
 there is an equivalence class $[\sigma]$ of section classes in $P_C$ and an element $b\in H_*(F)$  such that:
$$0\neq \sum_{\sigma'\in[\sigma]}\la \iota_C(a),\iota_C(b)\ra_{\sigma',2}^{P_C}=\sum_{\sigma'\in [\sigma]}\la [pt]; \iota_C(a),\iota_C(b),\iota_C([F]),...,\iota_C([F])\ra_{\sigma',l}^{P_C},$$
where the last equality follows from the Divisor axiom for Gromov-Witten invariants.
Applying the product formula \eqref{simpversionPF} as in the proof of Theorem \ref{unirulingfib2} we conclude that for every $a\in H_*(F)$ there is a non-vanishing Gromov-Witten invariant in $P$, namely:
\begin{equation}\label{preuveunireglage}
\la D; \iota(a),\iota(b),\pi^{-1}(c_3^B),...,\pi^{-1}(c^B_l)\ra_{\iota_{P_C}(\sigma),l}^P\neq 0.
\end{equation}
Taking $a=\ov{a}$ gives a contradiction, hence $\ker\iota$ is trivial. Now, the Corollary  follows directly from Corollary \ref{WeinFib}.
\qed\\ 

\bibliographystyle{plain}
\bibliography{bibarticle}

\begin{thebibliography}{10}

\bibitem{Bl}
Andr{\'e} Blanchard.
\newblock Sur les vari\'et\'es analytiques complexes.
\newblock {\em Ann. Sci. Ecole Norm. Sup. (3)}, 73:157--202, 1956.

\bibitem{CL}
Bohui Chen and An-Min Li.
\newblock Symplectic virtual localization of {G}romov-{W}itten invariants,
  2006.
\newblock Preprint, arXiv.org:math/0610370.

\bibitem{HV}
H.~Hofer and C.~Viterbo.
\newblock The {W}einstein conjecture in the presence of holomorphic spheres.
\newblock {\em Comm. Pure Appl. Math.}, 45(5):583--622, 1992.

\bibitem{HLR}
Jianxun Hu, Tian-Jun Li, and Yongbin Ruan.
\newblock Birational cobordism invariance of uniruled symplectic manifolds.
\newblock {\em Invent. Math.}, 172(2):231--275, 2008.

\bibitem{H}
C.~Hyvrier.
\newblock A product formula for gromov-witten invariants, 2009.
\newblock Preprint, arXiv.org:math/0904.1492, accepted for publication in the
  Journal of Symplectic Geometry.

\bibitem{Kedra}
Jaros{\l}aw K{\c{e}}dra.
\newblock Restrictions on symplectic fibrations.
\newblock {\em Differential Geom. Appl.}, 21(1):93--112, 2004.
\newblock With an appendix by the author and Kaoru Ono.

\bibitem{Kollar}
J{\'a}nos Koll{\'a}r.
\newblock {\em Rational curves on algebraic varieties}, volume~32 of {\em
  Results in Mathematics and Related Areas. 3rd Series. A Series of Modern
  Surveys in Mathematics}.
\newblock Springer-Verlag, Berlin, 1996.

\bibitem{Kollar1}
J{\'a}nos Koll{\'a}r.
\newblock Low degree polynomial equations: arithmetic, geometry and topology.
\newblock In {\em European {C}ongress of {M}athematics, {V}ol.\ {I}
  ({B}udapest, 1996)}, volume 168 of {\em Progr. Math.}, pages 255--288.
  Birkh\"auser, Basel, 1998.

\bibitem{LM}
F.~Lalonde and D.~McDuff.
\newblock Symplectic structures on fiber bundles.
\newblock {\em Topology}, 42(2):309--347, 2003.

\bibitem{LMP}
F.~Lalonde, D.~McDuff, and L.~Polterovich.
\newblock Topological rigidity of {H}amiltonian loops and quantum homology.
\newblock {\em Invent. Math.}, 135(2):369--385, 1999.

\bibitem{LT}
J.~Li and G.~Tian.
\newblock Virtual moduli cycles and {G}romov-{W}itten invariants of general
  symplectic manifolds.
\newblock In {\em Topics in symplectic $4$-manifolds (Irvine, CA, 1996)}, First
  Int. Press Lect. Ser., I, pages 47--83. Int. Press, Cambridge, MA, 1998.

\bibitem{RuanLi}
{T.-J.} {Li} and Y.~{Ruan}.
\newblock Uniruled symplectic divisors, 2007.
\newblock Preprint, arXiv.org:math/0711.4254.

\bibitem{LiuT}
Gang Liu and Gang Tian.
\newblock Weinstein conjecture and {GW}-invariants.
\newblock {\em Commun. Contemp. Math.}, 2(4):405--459, 2000.

\bibitem{Lu}
Guangcun Lu.
\newblock The {W}einstein conjecture in the uniruled manifolds.
\newblock {\em Math. Res. Lett.}, 7(4):383--387, 2000.

\bibitem{MS2}
D.~McDuff and D.~Salamon.
\newblock {\em Introduction to symplectic topology}.
\newblock Oxford Mathematical Monographs. The Clarendon Press Oxford University
  Press, New York, second edition, 1998.

\bibitem{MS}
D.~McDuff and D.~Salamon.
\newblock {\em {$J$}-holomorphic curves and symplectic topology}, volume~52 of
  {\em American Mathematical Society Colloquium Publications}.
\newblock American Mathematical Society, Providence, RI, 2004.

\bibitem{MDuniruled}
Dusa McDuff.
\newblock Hamiltonian {$S^1$}-manifolds are uniruled.
\newblock {\em Duke Math. J.}, 146(3):449--507, 2009.

\bibitem{Rvirt}
Y.~Ruan.
\newblock Virtual neighborhoods and pseudo-holomorphic curves.
\newblock In {\em Proceedings of 6th G\"okova Geometry-Topology Conference},
  volume~23, pages 161--231, 1999.

\bibitem{RT}
Y.~Ruan and G.~Tian.
\newblock A mathematical theory of quantum cohomology.
\newblock {\em J. Differential Geom.}, 42(2):259--367, 1995.

\bibitem{Se}
P.~Seidel.
\newblock {$\pi\sb 1$} of symplectic automorphism groups and invertibles in
  quantum homology rings.
\newblock {\em Geom. Funct. Anal.}, 7(6):1046--1095, 1997.

\bibitem{V}
Claude Viterbo.
\newblock A proof of {W}einstein's conjecture in {${\bf R}^{2n}$}.
\newblock {\em Ann. Inst. H. Poincar\'e Anal. Non Lin\'eaire}, 4(4):337--356,
  1987.

\bibitem{Wein}
Alan Weinstein.
\newblock On the hypotheses of {R}abinowitz' periodic orbit theorems.
\newblock {\em J. Differential Equations}, 33(3):353--358, 1979.

\end{thebibliography}

\end{document}